\documentclass[11pt,reqno]{amsart}
\usepackage{amsmath,amssymb,amsthm,upref,graphicx,mathrsfs}

\usepackage{color,xcolor}
\usepackage[
  colorlinks=true,
  linkcolor=blue,
  citecolor=blue,
  urlcolor=blue]{hyperref}

\newtheorem{thm}{Theorem}[section]

\newtheorem{lem}[thm]{Lemma}

\numberwithin{equation}{section}

\begin{document}

\leftline{ \scriptsize}

\vspace{1.3 cm}
\title
{Chollet's permanent conjecture for $4\times 4$ matrices}
\author{ Kijti Rodtes $^{\ast}$ }
\thanks{{\scriptsize
		\newline MSC(2010): 15A15, 15B48.  \\ Keywords: Chollet's permanent conjecture, Correlation matrices \\
		E-mail addresses:  kijtir@nu.ac.th (Kijti Rodtes).\\
		$^{*}$ Department of Mathematics, Faculty of Science, Naresuan University, Phitsanulok 65000, Thailand.\\
		\\}}
\hskip -0.4 true cm

\maketitle


\begin{abstract} In the year 1982, John Chollet conjectured that, for any pair of $n\times n$ positive semidefinite matrices $A,B$,  $per(A)\cdot per(B)\geq per(A\circ B)$, where $A\circ B$ is the Hardamard product of $A$ and $B$.  This conjecture was proved to be valid for $n=2, 3$ in the year 1987.  In this paper, we show that the conjecture holds true for $n=4$. 
\end{abstract}

\vskip 0.2 true cm


\pagestyle{myheadings}
\markboth{\rightline {\scriptsize Kijti Rodtes}}
{\leftline{\scriptsize }}
\bigskip
\bigskip


\vskip 0.4 true cm

\section{Introduction}

John Chollet proposed the conjecture analogue the Oppenheim's inequality concerning permanents for any pair of $n\times n$ positive semidefinite matrices $A,B$ in the year 1982, \cite{Chollet}, as:
$$  per(A)\cdot per(B)\geq per(A\circ B),$$ where $A\circ B$ is the Hardamard product of $A$ and $B$.  He also reduced the problem into the following form: for any $n\times n$ positive semidefinite matrix $A$,
\begin{equation} \label{conjec}
[per(A)]^2\geq per(A\circ \bar{A}),
\end{equation}
where $\bar{A}$ denotes the conjugate of $A$. It is also well known (see, e.g., \cite{GM} and \cite{Hut}), by using the standard scaling arguments, that it suffices to consider the conjecture (\ref{conjec}) only for the case of $n\times n$ correlation matrices.  Recall that, a correlation matrix $A=(A_{ij})$ is a positive semidefinite matrix for which each diagonal entry is $1$; which also implies that all $|A_{ij}|\leq 1$.

The conjecture was proved to be valid for $n=2, 3$ by R. J. Gregorac and I.R. Hentzel in the year 1987, \cite{GH}.  There is recently a progress to this conjecture for $4\times 4$ real matrices by G. Hutchinson in the year 2021, \cite{Hut}.    In this paper, we show that:
\begin{thm}
	The conjecture (\ref{conjec}) holds true for $n=4$.
\end{thm}  
In order to prove the result, we make use the following facts:
\begin{enumerate}
		\item Lieb inequality \cite{Lieb}: $per(A)\geq per(B)\cdot per(D)$, for any positive semidefinite matrix $A$ in the block form: $A=\begin{bmatrix}
	B&C \\
	C^*&D
	\end{bmatrix}$, where $B,D$ are square blocks.  Since $B,D$ are positive semidefinite,   $per(A)\geq per(B)\cdot per(D)\geq 0$. \\
	\item A result of Grone and Pierce (Theorem 3 in \cite{GP}): If $A=(A_{ij})$ is a correlation matrix of size $n\times n$, then $$per(A)\geq \frac{1}{n}\sum_{i=1}^n\sum_{j=1}^n|A_{ij}|^2. $$
	
\end{enumerate} 

\section{Proof}
Throughout the proof, we fix the notions of the following correlation matrix $A$ and its entries:  $$A:=\begin{pmatrix}
1 & x & y&z\\
\bar{x} & 1 & t&u \\
\bar{y} & \bar{t} & 1& w\\
\bar{z} & \bar{u} & \bar{w}&1 
\end{pmatrix},$$
where $x,y,z,t,u,w$ are complex numbers.  It is a direct calculation to see that
$$ per(A)=1+|x|^2+|y|^2+|z|^2+|t|^2+|u|^2+|w|^2+|xw|^2+|yu|^2+|zt|^2+\sum_{i=1}^7y_i, $$
and $$per(A\circ \bar{A})=1+|x|^4+|y|^4+|z|^4+|t|^4+|u|^4+|w|^4+|xw|^4+|yu|^4+|zt|^4+\sum_{i=1}^7z_i, $$
where, ($Re(c)$:=real part of the complex number $c$),  $$ \begin{matrix}
y_1:=2Re(x\bar{y}t), &y_2:=2Re(xu\bar{z}), &y_3:=2Re(yw\bar{z}),& y_4:=2Re(tw\bar{u}) \\
y_5:=2Re(xtw\bar{z}),& y_6:=2Re(y\bar{t}u\bar{z}), & y_7:=2Re(xu\bar{y}\bar{w}), &
\end{matrix}  $$
and 
$$ \begin{matrix}
z_1:=2|xyt|^2, &z_2:=2|xuz|^2, &z_3:=2|ywz|^2,& z_4:=2|twu|^2 \\
z_5:=2|xwzt|^2,& z_6:=2|yuzt|^2, & z_7:=2|xwyu|^2. &
\end{matrix}  $$
To prove the result, the following three lemmas are also needed.
\begin{lem}\label{key1} 
	By using the above notations,  
	\begin{eqnarray*}
	1+y_1&\geq& |x|^2+|y|^2+|t|^2:=s(y_1), \\
	1+y_2&\geq& |x|^2+|z|^2+|u|^2:=s(y_2), \\
	1+y_3&\geq& |y|^2+|z|^2+|w|^2:=s(y_3), \\
	1+y_4&\geq& |t|^2+|u|^2+|w|^2:=s(y_4) .
	\end{eqnarray*}
\end{lem}
\begin{proof}
For each $i=1,2,3,4$, denote $A(i)$ the $3\times 3$ principal submatrix of $A$ removing the $i$-th row and the $i$-th column.  Note that $A(i)$ is positive semidefinite and thus $\det(A(i))\geq 0$.  The inequality involving $1+y_i$ follows by evaluating $\det(A(5-i))\geq 0$ directly. 	
\end{proof}

\begin{lem}\label{key2} 
	By using the above notations, 
	\begin{eqnarray*}
		per(A)&\geq& 1+y_1+s(y_1)\geq 0, \\
		per(A)&\geq& 1+y_2+s(y_2)\geq 0, \\
		per(A)&\geq& 1+y_3+s(y_3)\geq 0, \\
		per(A)&\geq& 1+y_4+s(y_4)\geq 0 .
	\end{eqnarray*}
\end{lem}
\begin{proof}
By Lieb inequality, we have that $$per(A)\geq per(A(4))\cdot per([1])\geq 0 \;\;\hbox{ and }\;\; per(A)\geq per([1])\cdot per(A(1))\geq 0.$$  By direct calculation, $per(A(4))=1+y_1+s(y_1)$ and $per(A(1))=1+y_4+s(y_4)$, which yields the first row and the forth row inequalities, respectively. 	For the second row and the third row inequalities, we apply the Lieb inequality to the positive semidefinite matrices $A_2:=P_{(3\;4)}AP_{(3\;4)}$ and $A_3:=P_{(1\;2)}AP_{(1\;2)}$, respectively, where $P_{\sigma}$ is the permutation matrix associated to the permutation $\sigma=(1\; 2)$ (and $\sigma=(3\; 4)$).  Explicitly, $per(A)=per(A_2)=per(A_3)$, $per(A_2(4))=1+y_2+s(y_2)$, $per(A_3(1))=1+y_3+s(y_3)$  and
$$per(A_2)\geq per(A_2(4))\cdot per([1])\geq 0\;\; \hbox{ and }\;\; per(A_3)\geq per([1])\cdot per(A_3(1))\geq 0.$$  	
\end{proof}

\begin{lem}\label{key3} 
	For real numbers $a,b,c$, if $a^2+b^2+c^2\leq 1$, then
$$a^2+b^2+c^2\geq a^4+b^4+c^4+2(a^2b^2+a^2c^2+b^2c^2) .$$
\end{lem}
\begin{proof}
	The inequality follows by the fact that $x\geq x^2$ for any real number $x\in [0,1]$.
\end{proof}
We now divide the proof into 3 cases regarding the signs of $y_1,y_2,y_3,y_4$.\\

\textbf{Case 1: all $y_1,y_2,y_3,y_4$ are non-negative.} \\

Let $s(y_{max}):=\max\{s(y_1),s(y_2),s(y_3),s(y_4)\}$.
By Lemma \ref{key1} and Lemma \ref{key2}, 
\begin{eqnarray*}
[per(A)]^2 &\geq & [(1+y_{max})+s(y_{max})]^2 \\
 &=& (1+y_{max})^2+2(1+y_{max})s(y_{max})+(s(y_{max}))^2\\
 &\geq& (1+y_{max})^2+2(s(y_{max}))s(y_{max})+(s(y_{max}))^2\\
 &=& (1+y_{max})^2+\frac{1}{2}(4(s(y_{max}))^2)+(s(y_{max}))^2.
\end{eqnarray*}
Since $y_{max}\geq 0$ (assumption of Case 1), we have that $(1+y_{max})^2\geq 1$.  By the definition of $s(y_{max})$, it is clear that $(s(y_{max}))^2\geq (s(y_i))^2$ for $i=1,2,3,4$. Also, since $A$ is a correlation matrix, $|\alpha|^2\leq 1$ for all $\alpha\in \{x,y,z,t,u,w\}$.  This yields that $s(y_{max})\geq |xw|^2+|yu|^2+|zt|^2$.   So,
\begin{equation}\label{boundcase1}
[per(A)]^2\geq 1+\frac{1}{2}[(s(y_1))^2+(s(y_2))^2+(s(y_3))^2+(s(y_4))^2]+(|xw|^2+|yu|^2+|zt|^2)^2.
\end{equation}
Note that 
\begin{eqnarray*}
\frac{1}{2}(s(y_1))^2&=&\frac{1}{2}(|x|^4+|y|^4+|t|^4)+(|xy|^2+|xt|^2+|yt|^2)\\
&\geq& \frac{1}{2}(|x|^4+|y|^4+|t|^4)+(|xyt|^2+|xyt|^2+|xyt|^2)\\
&\geq& \frac{1}{2}(|x|^4+|y|^4+|t|^4)+z_1.
\end{eqnarray*}
Here, we have used $|t|^2,|y|^2,|x|^2\leq 1$ to conclude that $|xy|^2\geq |xyt|^2, |xt|^2\geq |xyt|^2$ and $|yt|^2\geq |xyt|^2$.
Similarly, 
\begin{eqnarray*}
\frac{1}{2}(s(y_2))^2&\geq& \frac{1}{2}(|x|^4+|z|^4+|u|^4)+z_2, \\
\frac{1}{2}(s(y_3))^2&\geq& \frac{1}{2}(|y|^4+|z|^4+|w|^4)+z_3, \\
\frac{1}{2}(s(y_4))^2&\geq& \frac{1}{2}(|t|^4+|u|^4+|w|^4)+z_4. 
\end{eqnarray*}
Moreover,
\begin{eqnarray*}
(|xw|^2+|yu|^2+|zt|^2)^2&=&(|xw|^4+|yu|^4+|zt|^4)+2(|xwyu|^2+|xwzt|^2+|yuzt|^2)\\
&=&(|xw|^4+|yu|^4+|zt|^4)+z_5+z_6+z_7.
\end{eqnarray*}  
Substituting these inequalities into (\ref{boundcase1}), we conclude that $[per(A)]^2\geq per(A\circ \bar{A})$.\\

\textbf{Case 2: all $y_1,y_2,y_3,y_4$ are negative.} \\

By a result of Grone and Pierce (Theorem 3 in \cite{GP}), we compute that
$$per(A)\geq 1+\frac{1}{2}(|x|^2+|y|^2+|z|^2+|t|^2+|u|^2+|w|^2). $$
By setting $T:=|x|^2+|y|^2+|z|^2+|t|^2+|u|^2+|w|^2$, the above inequality becomes 
\begin{equation}\label{bound2}
[per(A)]^2\geq  1+T+\frac{1}{4}T^2.
\end{equation}
By the direct expansion of $T^2$, ignoring unwanted terms (all are non-negative) and using the arithmetic-geometric means inequality, we compute that
\begin{eqnarray*}
		\frac{1}{4}T^2&\geq  &\frac{1}{4}( |x|^4+|y|^4+|z|^4+|t|^4+|u|^4+|w|^4)+\frac{1}{2}(|xw|^2+|yu|^2+|zt|^2) \\
		&=& \frac{1}{2}(\frac{1}{2} (|x|^4+|w|^4)+\frac{1}{2}(|z|^4+|t|^4)+\frac{1}{2}(|y|^4+|u|^4))+\frac{1}{2}(|xw|^2+|yu|^2+|zt|^2) \\
		&\geq& \frac{1}{2}(|xw|^2+|yu|^2+|zt|^2)+\frac{1}{2}(|xw|^2+|yu|^2+|zt|^2)\\
		&=& |xw|^2+|yu|^2+|zt|^2.
	\end{eqnarray*}
  Since $y_1<0$, by Lemma \ref{key1}, $s(y_1)<1$ and thus  $|xw|^2+|yu|^2+|zt|^2<1$ (in fact, any $y_i<0$ yields the same conclusion).  Now, by Lemma \ref{key3}, 
\begin{eqnarray*}
\frac{1}{4}T^2 &\geq&  |xw|^2+|yu|^2+|zt|^2  \\
&\geq& (|xw|^4+|yu|^4+|zt|^4)+2(|xwyu|^2+|xwzt|^2+|yuzt|^2)\\
&=&(|xw|^4+|yu|^4+|zt|^4)+z_5+z_6+z_7.  
\end{eqnarray*}
Also, by using Lemma \ref{key1} and the assumption that all $y_1,y_2,y_3,y_4$ are negative, it turns out that $s(y_i)<1$ for all $i=1,2,3,4$.  After applying Lemma \ref{key3} to each $s(y_i)$ and summing up the 4 inequalities, we obtain that
\begin{eqnarray*}
2T&\geq &2(|x|^4+|y|^4+|z|^4+|t|^4+|u|^4+|w|^4)+2(|xy|^2+|xt|^2+|yt|^2) \\
&+&2(|xz|^2+|xu|^2+|zu|^2)+2(|yz|^2+|yw|^2+|zw|^2)+2(|tu|^2+|tw|^2+|uw|^2) \\
&\geq& 2(|x|^4+|y|^4+|z|^4+|t|^4+|u|^4+|w|^4)+6(|xyt|^2+|xzu|^2+|yzw|^2+|twu|^2),
\end{eqnarray*}
(we have used $|\alpha|^2\leq 1$ for any $\alpha\in \{x,y,z,t,u,w\}$), which means that $$ T\geq |x|^4+|y|^4+|z|^4+|t|^4+|u|^4+|w|^4+\frac{3}{2}(z_1+z_2+z_3+z_4). $$
Substituting these inequalities into (\ref{bound2}), one gets that $[per(A)]^2\geq per(A\circ \bar{A})$.\\

\textbf{Case 3: $\{y_1,y_2,y_3,y_4\}$ contains both negative and non-negative real numbers.} \\

Let $s(y_m):=\max\{s(y_i)\;|\; y_i\geq 0\}$.  If $s(y_m)\geq 1$, then, by Lemma \ref{key1}, $1+y_m\geq s(y_m)$ and $s(y_m)\geq s(y_i)$ for each $i=1,2,3,4$ (if $y_i<0$, then $s(y_i)<1$).  We thus proceed the proof as in Case 1.  However, if  $s(y_m)<1$, then $s(y_i)<1$ for all $i=1,2,3,4$.  We thus  proceed the proof as in Case 2.

The proof is now completed.

\section*{Acknowledgments}
The author would like to thank anonymous referee(s) for reviewing this manuscript.  He also would like to thank Fundamental Research Fund (Naresuan University: FF 2566) under the project number FRB660001/0179 for financial support.

\end{document}